\documentclass[a4paper]{article}
\usepackage[T1]{fontenc}
\usepackage[english,francais]{babel}
\usepackage{amssymb,url,xspace,amsthm}
\usepackage{makeidx}
\usepackage[latin1]{inputenc}
\newtheorem{theor}{Théorème}
{\theoremstyle{definition} }
\newtheorem{pr}{Proposition}
\newtheorem{lemme}{Lemme}
{\theoremstyle{remark} \newtheorem{re}{Remarque}}
{\theoremstyle{definition} }
{\theoremstyle{definition} \newtheorem*{prc}{Proposition}}
{\theoremstyle{definition} \newtheorem*{thc}{Théorème}}
{\theoremstyle{definition} }
{\theoremstyle{remark} }
{\theoremstyle{definition} \newtheorem{de}{Définition}}

\providecommand{\bysame}{\leavevmode ---\ }
\providecommand{\og}{``}
\providecommand{\fg}{''}
\providecommand{\smfandname}{et}
\providecommand{\smfedsname}{\'eds.}
\providecommand{\smfedname}{\'ed.}
\providecommand{\smfmastersthesisname}{M\'emoire}
\providecommand{\smfphdthesisname}{Th\`ese}

\begin{document}
\selectlanguage{francais}
\author{Jérôme Depauw}
\title{Théorème ergodique pour cocycle harmonique, applications au milieu aléatoire\\
 Ergodic theorem for harmonic cocycle,  applications in random environment}
\date{18 janvier 2013}
\maketitle

\abstract{Dans ce travail est démontré le théorème ergodique ponctuel pour cocycle harmonique de degré~$1$
d'une action mesurable stationnaire de ${\mathbb Z}^d$ sur un espace de probabilité.
Ceci constitue un prolongement de l'article de Boivin et Derriennic (1991), qui portait sur les cocycles non nécessairement
harmoniques. L'hypothèse d'harmonicité permet, dans le cas elliptique, d'abaisser la condition d'intégrabilité à $L^2$,
alors que Boivin et Derriennic ont montré que la condition 
optimale dans le cas non harmonique est la finitude de la norme de Lorentz $L_{d,1}$. Ils ont montré notamment 
que $L^d$ ne suffit pas. 
Le présent travail constitue aussi une suite à un article de Berger et Biskup (2007) qui portait
sur le cas harmonique non elliptique, en dimension $d=2$.
Enfin des applications de ce théorème au milieu aléatoire sont présentées. 

\textbf{Abstract}. In this work we prove the pointwise ergodic theorem for harmonic degree 1 cocycle of a measurable stationary 
action of $ {\mathbb Z} ^ d $ on a probability space.
In a precedent  paper  Boivin and Derriennic  (1991) studied  this theorem for not necessarily harmonic cocycles. The harmonic hypothesis allows, in the elliptic case, to change the integrability condition to $L^2$, while Boivin and Derriennic showed that the optimum condition in the non-harmonic case is the finiteness of Lorentz's norm $L_{d,1}$. They showed in particular that $L^d$ is not enough. Berger and Biskup published in 2007 a paper on the harmonic 
 not elliptic case, but only in dimension $d=2$. Finally, applications of this theorem in random media are presented.}


\section{Introduction}
L'objet de ce travail est le théorème ergodique ponctuel pour
cocycle harmonique d'une action stationnaire de ${\mathbb Z}^d$ (théorème~\ref{thharmo}
ci-dessous). La présentation du cadre dans lequel se place ce résultat fait l'objet de cette introduction.
C'est aussi l'occasion de présenter deux applications
de ce théorème.
Cette introduction est divisée en trois paragraphes. 
Le premier paragraphe propose quelques
rappels sur le théorème ergodique en dimension~$1$. Le second porte  sur
la notion de  cocycle de degré~$1$ en dimension~$d$ et sur le théorème
ergodique associé.
Enfin dans le troisième nous présentons la notion de cocycle harmonique, le théorème ergodique qui fait
l'objet de cet article, et son rôle dans certaines questions issues du domaine du milieu aléatoire.
\subsection{Théorème ergodique en dimension $d=1$}
Soit un espace de probabilité $(\Omega,{\mathcal A},P)$, et $T$ une transformation mesurable, inversible,
telle que $T$ et $T^{-1}$
préservent la probabilité $P$. Le théorème ergodique ponctuel  de Birkhoff énonce que pour toute fonction intégrable $f$
sur $\Omega$, les moyennes $\frac{1}{n}\sum_{k=0}^{n-1}f\circ T^k$ convergent presque sûrement,
lorsque $n$ tend vers l'infini. 
Si  la transformation $T$ est ergodique,
c'est-à-dire n'admet pas d'ensemble invariant de probabilité non égale à $0$ ou $1$, la limite  est constante, égale à l'intégrale $\int_\Omega f \ dP$ de la fonction $f$.

Ce résultat, ainsi que ses différentes démonstrations, sont devenues très classiques. 
 Il est cependant utile, afin de rendre plus naturel le passage à la  dimension supérieure, 
de rappeler une de ces démonstrations. Celle présentée ici pour $T$ ergodique consiste en trois étapes.

 On constate d'abord
 que la convergence est aisée à vérifier
sur les cobords, c'est-à-dire sur les fonctions $f$ pouvant s'écrire $f=g\circ T-g$, avec $g\in L^1$.
En effet dans ce cas, les termes de la somme $S_n(\omega)=\sum_{k=0}^{n-1}f\circ T^k(\omega)$   se simplifiant deux à deux, celle-ci s'écrit $S_n=g\circ T^n-g$. 
La convergence ponctuelle des moyennes de Césaro $\frac{1}{n}S_n$ vers $0$,
 qui se résume donc à celle du terme $\frac{1}{n}g\circ T^n$, est alors une application directe du théorème de Borel-Cantelli.
 En effet celui-ci assure qu'une condition suffisante à cette converge presque sûre est que la série
 $$\sum_{n=1}^{+\infty}P\Bigl(\Bigl|\frac{1}{n}g\circ T^n\Bigr|>\varepsilon\Bigr)
 $$
 converge pour tout $\varepsilon$. Par invariance de $P$ sous l'action de $T$, cela revient à la convergence de la série
$$\sum_{n=1}^{+\infty}P(|g|>n\varepsilon).
 $$
Ceci est bien une condition équivalente à la convergence de l'intégrale $\int_\Omega |g| \ dP$.

La deuxième étape de la démonstration du théorème ergodique
consiste à vérifier que toute 
fonction intégrable de moyenne nulle peut être approchée dans $L^1$ par des cobords. Cette propriété
de densité peut être déduite du théorème de Hann-Banach. En effet d'après ce théorème, 
 s'il existait une fonction $f$ d'intégrale nulle qui ne soit pas dans l'adhérence de l'espace des cobords,
alors il existerait une forme linéaire nulle sur l'espace des cobords et non nulle sur $f$.
Une forme linéaire de $L^1$ s'identifiant  à une fonction $h\in L^{\infty}$,
on aurait donc
$$   \int (g\circ T-g)h\ dP=0\mbox{ pour tout }g\in L^1,\qquad\mbox{et }\int fh\ dP\neq 0. 
$$
De la première condition et de l'ergodicité de $T$ on déduit que $h$ est constante,
la deuxième condition contredit alors le fait que $f$ est d'intégrale non nulle. Ceci prouve
la densité cherchée.

La  démonstration du théorème ergodique ponctuelle s'achève en montrant que l'ensemble
des fonctions $f$ vérifiant le théorème ergodique ponctuel
est fermé dans $L^1$. Ceci est une conséquence de l'inégalité maximale
$$P\Bigl(\sup_n\Bigl|\frac{1}{n}S_n\Bigr|>\varepsilon \Bigr)<\frac{\int_\Omega |f|\ dP}{\varepsilon}
$$
dont  nous ne détaillons pas la démonstration ici. Cette inégalité maximale est  le point délicat de la démonstration du théorème ergodique
ponctuel.
\subsection{Théorème ergodique en dimension $d\geq 1$}
Il y a plusieurs généralisations possibles du théorème ergodique au cas multidimensionnel. 
La plus connue est celle de Wiener, qui consiste à faire des moyennes sur des cubes de dimension~$d$ dans ${\mathbb Z}^d$ (voir~\cite{Wiener}).


Une autre généralisation du théorème (et de
la démonstration) exposé ci-dessus à la dimension $d\geq 1$  est la  suivante.
Soient les  transformations  $T_1,\ldots , T_d$ mesurables inversibles, 
préservant la probabilité $P$. Lorsque de plus elles commutent deux à deux, on obtient une 
action de ${\mathbb Z}^d$ en considérant  pour $\vec n=(n_1,\ldots,n_d)\in{\mathbb Z}^d$
la transformation $T_{\vec n}=T_1^{n_1}\circ\cdots\circ T_d^{n_d}$. Commençons, par analogie
avec la démonstration en dimension~$1$, à nous intéresser
à la convergence ponctuelle des quantités  $\frac{1}{|\vec n|}S_{\vec n}(\omega)$ vers $0$, pour $|\vec n|$ tendant
vers l'infini,  où
\begin{itemize}  
\item $S_{\vec n}=g\circ T_{\vec n}-g$;
\item $|\vec n|=\sum_{i=1}^d|n_i|$.
\end{itemize}  
Celle-ci se résume, comme ci-dessus, à la convergence vers zéro du terme $\frac{1}{|\vec n|}g\circ T_{\vec n}$. 
 Il découle à nouveau du théorème de Borel-Cantelli qu'une condition suffisante à cette convergence est que la série
 $$\sum_{n_1=-\infty}^{+\infty}\cdots \sum_{n_d=-\infty}^{+\infty}P\Bigl(\Bigl|\frac{1}{|\vec n|}g\circ T_{\vec n}\Bigr|>\varepsilon\Bigr)
 $$
 converge pour tout $\varepsilon$. 
 Or pour tout entier $n$, le nombre de vecteurs $\vec n$ de norme $|\vec n|$ fixée, égale à $n$, est de l'ordre de $n^{d-1}$.
 Par invariance de $T$, la convergence de la série ci dessus  revient donc à la convergence de la série
$$\sum_{n=1}^{+\infty}n^{d-1}P(|g|>n\varepsilon),
 $$
ce qui est cette fois une condition équivalente à la convergence de l'intégrale $\int_\Omega |g|^d \ dP$.
Ainsi on voit que la condition d'intégrabilité pour la convergence ponctuelle étudiée dans ce paragraphe
dépend de la dimension de l'action: c'est $g\in L^d$.

Pour identifier l'adhérence de l'espace des sommes $S_{\vec n}$ considérées
ci-dessus, on remarque que celles-ci vérifient 
l'équation 
\begin{equation}\label{cc1}
S_{\vec n+\vec m}=S_{\vec n}+S_{\vec m}\circ T_{\vec n}.
\end{equation}
Nous introduisons donc la notion correspondante, par la définition suivante.
\begin{de}[Cocycle de degré $1$]
Un cocycle de degré 1 d'une action $T$  mesurable stationnaire de ${\mathbb Z}^d$
est un processus $(S_{\vec n})_{{\vec n}\in {\mathbb Z}^d}$ vérifiant
l'équation~(\ref{cc1}) ci-dessus.
\end{de}
Cette notion de cocycle de degré~$1$ d'une action multidimensionnelle, introduite ici de manière qui peut sembler détournée,
à l'occasion d'une démonstration, apparait en fait dans de nombreuses questions de
théorie ergodique. Citons les principales, sans vouloir être exhaustif.
\begin{enumerate}
\item Dans l'étude des systèmes dynamiques sous l'aspect de la cohomologie des groupes,
 les cocycles de degré $1\leq k\leq d$ ont été notamment
étudiés par Mackey~\cite{Mackey}, Feldman et Moore~\cite{FeldmanMoore1} et~\cite{FeldmanMoore2}, Katok et Katok~\cite{KatokKatok}, Depauw~\cite{DepauwETDS}.
\item Les cocycles de degré 1 d'action stationnaire de ${\mathbb Z}^d$ et à valeurs
dans  ${\mathbb R}^d$ joue un rôle particulier, car ils représentent une déformation aléatoire du réseau
${\mathbb Z^d}$, lorsque l'on  considère ce dernier plongé dans ${\mathbb R}^d$.
Ils interviennent dans la construction du flot spécial en dimension $d$ (voir~\cite{Katok}). 
\item\label{p3}
Ils apparaissent
aussi dans certaines démonstrations du théorème limite central  multidimensionnel pour marche aléatoire en milieu aléatoire.
En effet la déformation du réseau associée  permet de transformer ce type de marche aléatoire
en martingale. Or on sait
que la condition de martingale joue un rôle crucial dans le théorème limite
central. La littérature sur ce sujet est très abondante. Sur le rôle du théorème ergodique 
ponctuel pour cocycle de degré $1$ dans cette méthode de démonstration
du théorème limite central  on renvoie au survey~\cite{Biskup}, et notamment 
au paragraphe intitulé \og Sublinearity of the corrector\fg.
\item\label{p2} Enfin l'équation de cocycle de degré 1 a une interprétation physique
qui fait qu'elle apparait dans le domaine du milieu aléatoire, et notamment
des réseaux de conductances stationnaires. En effet elle est vérifiée par la différence de potentiel 
entre l'origine et le point $\vec n$, lorsque l'accroissement de ce potentiel le long des arêtes du réseau ${\mathbb Z}^d$
est stationnaire. Ce potentiel est  un moyen d'étudier la résistivité équivalente du milieu infini. Citons notamment  Golden et Papanicolaou~\cite{GoldenPapanicolaou}.
\end{enumerate}
Les points~\ref{p3} et~\ref{p2} sont exposés de manière plus détaillée dans le paragraphe suivant.

Ne considérant ici essentiellement que des cocycles de degré 1,  la spécification du degré sera souvent
sous entendue.

Revenons au problème du théorème ergodique ponctuel. Ce théorème est énoncé
ci-dessous (théorème de Boivin-Derriennic), mais commençons par rappeler, sans entrer dans
les détails, le principe de sa démonstration, 
qui suit le shéma vu dans le paragraphe précédent sur la dimension~$1$. 
Supposons que l'action $T$ est ergodique, c'est-à-dire qu'il n'y a pas d'ensemble invariant
par les $d$ transformations $T_1,\ldots T_d$ simultanément, sauf de mesure $0$ ou $1$.
On peut vérifier que les cobords, c'est-à-dire les
cocycles
s'écrivant sous la forme $S_{\vec n}=g\circ T_{\vec n}-g$, sont denses dans l'espace des cocycles
 d'intégrale nulle. Ensuite une inégalité maximale permet de montrer que les
 cocycles vérifiant le théorème ergodique forment un espace fermé pour la norme adéquate. Là encore,
il y a  une différence notable avec la dimension $d=1$. La norme majorante dans l'inégalité maximale n'est pas la norme $L^1$, ni même la norme $L^d$. Boivin et Derriennic ont montré dans~\cite{BoivinDerriennic} que 
la norme optimale est la norme de 
 Lorentz $L_{d,1}$, qui s'intercale entre les normes $L^d$ et $L^{d+\varepsilon}$.
 
Ces questions ne peuvent être détaillées ici.
Pour fixer les idées, donnons seulement un énoncé avec une condition d'intégrabilité plus forte que nécessaire, en terme
d'espace $L^p$. On aura besoin de considérer la base canonique de ${\mathbb R}^d$, notée
$(\vec e_1,\ldots,\vec e_d)$.
\begin{thc}[Boivin \& Derriennic]\label{thBD} Soit $(S_{\vec n})_{{\vec n}\in {\mathbb Z}^d}$ un cocycle d'une action $T$ mesurable stationnaire
ergodique 
de ${\mathbb Z}^d$. Si il existe $\varepsilon>0$ tel que les fonctions  $f_i=S_{\vec e_i}$
soient de norme $L^{d+\varepsilon}$ finie pour $i=1,\ldots,d$
alors on a la convergence ponctuelle
$$\frac{S_{\vec n}-\sum_{i=1}^dn_i(\int_\Omega f_i \ dP)}{|\vec n|}\longrightarrow 0
$$
pour $|\vec n|$ tendant vers l'infini. 
\end{thc}
\begin{re}\label{re1}
Les mêmes auteurs ont aussi montré que dès que les fonctions $f_i$ sont  intégrables, la convergence ci-dessus
a  lieu au sens de la norme  $L^1$.
\end{re}
\begin{re}\label{re2}
Le terme soustrait au numérateur est l'intégrale de la fonction $S_{\vec n}$. En effet, de l'équation de cocycle~(\ref{cc1})
et de la stationnarité de l'action $T$ on déduit
$$
\int_\Omega S_{\vec n+\vec m}\ dP=\int_\Omega S_{\vec n}\ dP+\int_\Omega S_{\vec m}\ dP.
$$
Il en découle immédiatement, en écrivant $\vec n=\sum_{i=1}^dn_i\vec e_i$, que
$$\int_\Omega S_{\vec n}\ dP=\sum_{i=1}^dn_i\int_\Omega f_i\ dP.
$$
\end{re}
%
\subsection{Théorème ergodique pour cocycles har\-mo\-ni\-ques}
Les cocycles évoqués dans les points~\ref{p3} et~\ref{p2} vérifient une propriété supplémentaire d'harmonicité.
Détaillons cela en commençant par le point~\ref{p3}. Soient $(\Omega,{\mathcal A},P)$ un espace de probabilité
et $T$ une action mesurable stationnaire ergodique
de ${\mathbb Z}^d$
engendrée par $d$ transformations commutantes $T_1,\ldots,T_d$. Soient $c_1,\ldots,c_d$ des fonctions mesurables $>0$. 
Pour chaque $\omega\in\Omega$, on considère la marche aléatoire sur ${\mathbb Z}^d$
consistant à sauter d'un point $\vec n$ à ses $2d$ plus proches voisins $\vec n\pm \vec e_i$, avec des probabilités
données par les fonctions $(c_i)_i$ par
\begin{eqnarray}
P_{\omega}(X_{k+1}=\vec n+\vec e_i\ /\ X_k=\vec n)&=&\frac{c_i(T_{\vec n}\omega)}{\bar c(T_{\vec n}\omega)};
\cr
P_{\omega}(X_{k+1}=\vec n-\vec e_i\ /\ X_k=\vec n)&=&\frac{c_i(T_{\vec n-\vec e_i}\omega)}{\bar c(T_{\vec n}\omega)},
\end{eqnarray}
où $\bar c$ est la normalisation appropriée: $\bar c=\sum_{i=1}^d(c_i+c_i\circ T_i^{-1})$. 

L'intérêt de définir une  marche aléatoire à partir des fonctions $c_i$  est que cette marche est
 l'analogue, en temps et en espace discret, d'une diffusion dans un milieu 
aléatoire déterminé par $\omega$. Cela se voit en écrivant 
le générateur ${\mathcal L}_\omega$ associé,
défini sur les fonctions $u$ de ${\mathbb Z}^d$ par l'espérance conditionnelle
$$
({\mathcal L}_\omega u)(\vec n)=E_\omega\bigl(u(X_{k+1})-u(X_{k})/X_k=\vec n\bigr),
$$
où l'indice $\omega$ est là pour rappeler que les probabilités de sauts définissant la marche
aléatoire $(X_k)_k$ dépendent de l'environnement aléatoire $\omega$.
Pour faire apparaitre ce lien avec les diffusions, il faut donner l'expression de ${\mathcal L}_\omega$
à l'aide des fonctions $c_i$. Notons
pour cela  $\tau_i$ la translation de ${\mathbb Z}^d$ de vecteur $\vec e_i$,
et  posons 
 $\partial_i=\tau_i-I$ et $\partial_i^*=\tau_i^{-1}-I$. Enfin, à toute fonction $c$ définie sur $\Omega$,
 et à tout $\omega\in\Omega$, associons la fonction $c_\omega$  définie sur ${\mathbb Z}^d$ par  
 $c_\omega(\vec n)=c(T_{\vec n}\omega)$. 
 Alors un calcul élémentaire permet de vérifier que, pour toute fonction $u$ définie sur ${\mathbb Z}^d$,
 on a
 \begin{equation}\label{ccharmo}
 {\mathcal L}_\omega u=\frac{1}{\bar c_\omega}\sum_{i=1}^d\partial^*_i\bigl(
 (c_i)_\omega(\partial_iu)\bigr).
 \end{equation}
 On perçoit sans peine l'analogie avec l'opérateur de diffusion de la chaleur associé à un milieu
 de conductivité thermique  $A$, et de capacité
 thermique  $\lambda$, usuellement écrit
 $${\mathcal L}:u\mapsto \frac{1}{\lambda}\mbox{\rm div}\bigl(A(\mbox{\rm grad}u)\bigr).$$
 La conductivité thermique de l'arête  $[\vec n, \vec n+\vec e_i]$, à environnement $\omega$ fixé,
  est pour nous égale à  
 $(c_i)_\omega(\vec n)=c_i(T_{\vec n}\omega)$, pour $i=1,\ldots , d$, 
 et la capacité
 thermique du n\oe ud $\vec n$ à $\bar c_\omega(\vec n)=\bar c(T_{\vec n}\omega)$.

Dans l'idée d'aborder la question du théorème limite central pour la marche aléatoire $(X_k)_k$,
il est naturel d'étudier les fonctions $u$ harmoniques, c'est-à-dire vérifiant ${\mathcal L}_\omega u=0$.
En effet, dans ce cas la marche aléatoire $Y$ définie sur ${\mathbb R}$ par $Y_k=u(X_k)$
est une martingale, et est donc bien placée pour vérifier le théorème limite central.
En fait, pour ne pas perdre la structure aléatoire stationnaire définissant les probabilités
de saut de la marche, on  cherche $u$ sous la forme d'un cocycle. Plus précisément
on cherche les cocycles $S_{\vec n}$ tels que pour tout $\omega$, la fonction 
$u(\vec n)=S_{\vec n}(\omega)$ annule l'opérateur ${\mathcal L}_\omega$.  En remplaçant dans~(\ref{ccharmo})
on obtient la notion suivante.
\begin{de}[Cocycle harmonique] Le cocycle $(S_{\vec n})$ de degré $1$ est harmonique s'il vérifie
$$\sum_{i=1}^d (T_i^{-1}-I)(c_if_i)=0,
$$
où on a noté $f_i=S_{\vec e_i}$.
\end{de}
Dans l'optique de démontrer le théorème limite central pour $(X_k)_k$, l'intérêt d'un théorème ergodique ponctuel  pour ce type de cocycle peut être expliqué rapidement de la façon
suivante. Admettons qu'un tel cocycle existe, avec de plus les conditions d'intégrales $\int_\Omega f_1\ dP=1$
et $\int_\Omega f_i\ dP=0$, pour $ 2\leq i\leq d$. Supposons de plus que ce cocycle
vérifie le théorème ergodique énoncé ci-dessus (théorème de Boivin-Derriennic). Considérons la martingale $Y$
définie par  $Y_k=S_{X_k}(\omega)$. Celle-ci vérifie
$$\frac{Y_k-X_k^{(1)}}{|X_k|}\to 0   
$$
pour $|X_k|$ tendant vers l'infini, 
où $X_k^{(1)}$ désigne la première coordonnée de $X_k$. Autrement dit, l'écart entre la martingale $Y_k$
et la coordonnée $X_k^{(1)}$ est négligeable devant $|X_k|$. On peut en déduire que si 
$Y_k$ vérifie le théorème limite central, il en est de même de $X_k^{(1)}$. 
Pour démontrer que $Y_k$ vérifie le théorème limite central, on peut utiliser un théorème limite central
pour martingale, comme le théorème de Brown~\cite{Brown}. Cette méthode 
de démonstration du théorème limite central pour $(X_k^{(1)})_k$, dite \og méthode de la  martingale\fg\
est détaillée dans~\cite{Kozlov}). 

En fait l'existence d'un cocycle harmonique d'intégrales données est prouvée 
dans le cas elliptique, c'est-à-dire le cas où les conductances sont bornées inférieurement et supérieurement par
des constantes strictement positives. Enonçons
ce résultat dû à  Kunnemann (voir~\cite{Kunnemann}).
\begin{thc}[Kunnemann]\label{thKo}
Supposons qu'il existe des constantes $0<a<b$ telles que $a<c_i<b$ pour $i=1,\ldots,d$.
Alors pour tout vecteur $\vec y=(y_1,\ldots,y_d)\in{\mathbb R}^d$ il existe un unique cocycle harmonique $(S_{\vec n})_{\vec n\in{\mathbb Z}^d}$ tel que les fonctions 
définies par $f_i=S_{\vec e_i}$
vérifient  $f_i\in L^2$ et $\int_\Omega f_i\ dP=y_i$, pour $ 1\leq i\leq d$.
\end{thc}
%
Notons que l'intégrabilité $L^2$ donné par ce théorème ne suffit pas pour que le théorème de Boivin-Derriennic
s'applique. Cette intégrabilité a été améliorée par Boivin pour la dimension $d=2$ (voir~\cite{Boivin})
 puis par l'auteur de cet
article pour $d\geq 3$ (voir~\cite{DepauwIHP}) de
la manière suivante.
 \begin{thc}[Boivin, Depauw]
Supposons qu'il existe des cons\-tan\-tes $0<a<b$ telles que $a<c_i<b$ pour $i=1,\ldots,d$.
Alors il existe un $p>2$ tel que le cocycle
donné par le théorème de Kunnemann ci-dessus vérifie   $f_i\in L^p$ pour $ 1\leq i\leq d$.
\end{thc}
Ce résultat de Boivin en dimension $d=2$ rend opérationnelle  la méthode de la martingale expliquée
 ci-dessus dans 
cette dimension. En effet c'est  l'intégrabilité qui permet d'appliquer le théorème de Boivin-Derriennic
 (voir~\cite{Boivin}).

L'objet de ce travail est de montrer que le cocycle de Kunnemann vérifie le théorème ergodique ponctuel en toute
dimension.
\begin{theor}\label{thharmo}
Supposons qu'il existe des constantes $0<a<b$ telles que $a<c_i<b$ pour $i=1,\ldots,d$.
Alors le cocycle
$(S_{\vec n})_{\vec n\in{\mathbb Z}^d}$ donné par le théorème de Kunnemann ci-dessus vérifie   
$$\frac{S_{\vec n}-\sum_{i=1}^dn_i(\int_\Omega f_i \ dP)}{|\vec n|}\longrightarrow 0
$$
presque surement, pour $|\vec n|$ tendant vers l'infini. 
\end{theor}
Ainsi la méthode de la martingale proposée par Kozlov pour démontrer le théorème limite  central pour ce type
de marche aléatoire fonctionne en toute dimension (notons que ce théorème limite  central est déjà démontré,
par Boivin et Depauw dans~\cite{BoivinDepauw}, mais avec une méthode spectrale plus compliquée).
 
Pour achever cette introduction, présentons rapidement l'autre application de ce théorème,
évoquée point~\ref{p2} page~\pageref{p2} ci-dessus. 
Considérons un réseau de conductances électriques aléatoires, tel que, à $\omega$
fixé, la conductance
de l'arête $[\vec n,\vec n+\vec e_i]$ soit $c_i(T_{\vec n}\omega)$.
Cherchons un cocycle $(S_{\vec n})_{\vec n}$ tel que   la fonction  $u$ définie sur le réseau ${\mathbb Z}^d$
par $u(\vec n)=S_{\vec n}(\omega)$ soit le potentiel électrique dans le réseau. Cela signifie
que la quantité $f_i(\omega)=S_{\vec e_i}(\omega)$
est l'accroissement de potentiel dans l'arête $[0,\vec e_i]$. Il s'ensuit par la loi de Ohm que
$c_i(\omega)f_i(\omega)$  est l'intensité de courant dans cette arête. L'équation
d'harmonicité du cocycle $S_{\vec n}$ de la définition~\ref{ccharmo} exprime alors la loi des n\oe uds, selon laquelle
la somme des intensités sur les $2d$ arêtes arrivant à un  n\oe ud donné est nulle.
Ainsi un cocycle harmonique défini un potentiel électrique, à accroissement stationnaire. 
Le théorème~\ref{thharmo} s'interprète donc comme montrant l'existence d'un potentiel à l'infini,
dans un réseau stationnaire elliptique de conductance. Ce résultat complète un de nos
articles précédents sur l'existence d'un flux de courant à l'infini. Ce dernier reposait sur
le théorème ergodique ponctuel pour cocycle de degré $d-1$. 
Ces deux résultats (existence du potentiel et du courant stationnaire, convergeant à l'infini
dans presque tous
les environnements)
 montre que la conductivité équivalente du milieu infini est visible dans presque toutes les réalisations de ce milieu.
On renvoie à~\cite{DepauwIHP}, ainsi qu'à~\cite{DepauwMP} pour la version pour milieu continu.

\section{Démonstration du théorème~\ref{thharmo}}
La démonstration du théorème~\ref{thharmo} proposée ci-dessous ne suit pas un raisonnement par densité des
cobords comme
celles présentées dans l'introduction. 
On se base d'une part sur le théorème ergodique en dimension~$1$,
qui, dans le cadre de la dimension~$d$, assure la convergence cherchée
 dans toutes les directions rationnelles, et d'autre
part sur la régularité de Hölder des fonctions harmoniques  sur les réseaux. Pour cette régularité, 
qui est une version discrète du célèbre résultat de de Giorgi
(voir~\cite{DeGiorgi}), nous nous appuierons
sur un l'article de Delmotte~\cite{Delmotte}.

 Cette continuité höldérienne a pour conséquence la proposition suivant.
\begin{pr}\label{prharmo}
Supposons qu'il existe des constantes $0<a<b$ telles que $a<c_i<b$ pour $i=1,\ldots,d$.
Alors il existe $\alpha>0$ vérifiant:
pour tout  cocycle harmonique  $(S_{\vec n})_{\vec n}$ tel que les fonctions $f_i=S_{\vec e_i}$ soient
dans $L^2$ il existe, avec probabilité~$1$, une constante $C$
telle que pour tout $R\geq 1$ et tout $\vec m$, $\vec n\in{\mathbb Z}^d$ on a
$$\mbox{si }|\vec m|<R\mbox{ et }|\vec n|<R \mbox{ alors }|S_{\vec m}(\omega)-S_{\vec n}(\omega)|\leq RC
\Bigl( \frac{|\vec m-\vec n|}{R}\Bigr)^\alpha.
$$
\end{pr}
\begin{re} La constante $C$, explicitée ligne~(\ref{CC}) ci-dessous, dépend de $\omega$.
Mais pour ne pas alourdir les notations, ce dernier est omis.
\end{re}
\par\noindent{\it Démonstration de la proposition~\ref{prharmo}. --- } Commençons par le lemme suivant.
On rappelle que $\partial_i$ et $\partial_i^*$ désignent les opérateurs aux différences finies présentés ci-dessus
lors du calcul de l'opérateur ${\mathcal L}_\omega$,
et que $|\cdot|$ est la norme introduite précédemment, égale à la somme des valeurs absolues des coordonnées.
Remarquons que pour deux points $\vec m,\vec n\in{\mathbb Z}^d$, l'entier $|\vec m-\vec n|$ est le nombre d'arêtes
nécessaire pour relier ces deux points. Notamment ces points sont voisins ssi $|\vec m-\vec n|=1$.
 Enfin pour tout entier $R\geq 1$, l'expression \og boule de rayon R\fg\ désigne
les points $\vec n$ vérifiant $|\vec n|\leq R$.
\begin{lemme}[Inégalité de Poincaré]\label{poincarre}  
Pour toute fonction  $u$  définie sur le réseau ${\mathbb Z}^d$
et tout entier $R\geq 1$ on a 
$$\sum_{|\vec n|\leq R} |u(\vec n)-\bar u_R|^2\leq R^2 4\sum_{|\vec n|\leq R+1} \sum_{i=1}^d\bigl(
|\partial_iu(\vec n)|^2
+|\partial_i^*u(\vec n)|^2\bigr)
,
$$  
où  $\bar u_R$ désigne la moyenne de $u$ sur la boule de rayon $R$.
\end{lemme}
\par\noindent{\it Démonstration de lemme. --- } Ce résultat est des plus classiques. Rappelons en néanmoins la démonstration pour
la clarté de l'exposé. La boule  de rayon $R$ dans ${\mathbb Z}^d$ sera notée $B$ dans cette démonstration, et son
cardinal $|B|$. La somme à majorer s'écrit
$$
\sum_{|\vec n|\leq R} \Bigl|u(\vec n)-\frac{1}{|B|}\sum_{|\vec m|\leq R} u(\vec m)\Bigr|^2=
\sum_{|\vec n|\leq R} \Bigl|\frac{1}{|B|}\sum_{|\vec m|\leq R}\bigl(u(\vec n)- u(\vec m)\bigr)\Bigr|^2,
$$
ce qui, en posant le changement d'indice $\vec m\to\vec \ell$ défini par $\vec m =\vec n + \vec \ell$, est encore égale à 
$$
\sum_{|\vec n|\leq R} \Bigl|\frac{1}{|B|}\sum_{\vec \ell\atop |\vec n+\vec \ell|\leq R}\bigl(u(\vec n)- u(\vec n+\vec\ell)\bigr)\Bigr|^2.
$$
D'après l'inégalité de Cauchy-Schwarz appliquée à la somme intérieure, cela est majoré par
$$
\sum_{|\vec n|\leq R}\frac{1}{|B|}\sum_{\vec \ell\atop |\vec n+\vec \ell|\leq R}\bigl(u(\vec n)- u(\vec n+\vec\ell)\bigr)^2.
$$
Or on a
$$
u(\vec n)- u(\vec n+\vec\ell)=\sum_{\vec k\in\gamma_{\vec\ell}}\partial_{(i)}^{(*)}u(\vec n+\vec k),
$$
où $\gamma_{\vec \ell}$ désigne un chemin de $0$ à $\vec \ell$, c'est à dire une suite de sommets voisins
dont le premier est $0$ et le dernier est un voisin de $\vec n$; l'indice $i$, ainsi que 
la présence ou non de l'étoile, varient en fonction de la direction
de l'arête joignant les sommets voisins (c'est pourquoi ils sont notés entre parenthèses). On peut de plus choisir
ce chemin pour que  le nombre d'arêtes qui le composent soit égale à $|\vec \ell|$, et pour que les chemins $\vec n+\gamma_{\vec \ell}$, avec $\vec n$ et $\vec n+\vec\ell\in B$,   
soient inclus dans la boule  de rayon $R+1$,  ce que l'on suppose
fait dans la suite. On a donc, à nouveau d'après l'inégalité de Cauchy-Schwarz,
$$
\bigl(u(\vec n)- u(\vec n+\vec\ell)\bigr)^2\leq|\vec \ell|\sum_{\vec k\in\gamma_{\vec\ell}}\bigl(\partial_{(i)}^{(*)}u(\vec n+\vec k)\bigr)^2,
$$
ce qui peut être encore majoré par
$$|\vec \ell|\sum_{\vec k\in\gamma_{\vec\ell}}\sum_{i=1}^d\Bigl(\bigl(\partial_{i}^{*}u(\vec n+\vec k)\bigr)^2+\bigl(\partial_{i}u(\vec n+\vec k)\bigr)^2\Bigr).
$$
En divisant par $|B|$, en sommant sur $\vec n$, et en échangeant les sommes en $\vec k$ et $\vec n$, on a
\begin{eqnarray*}\sum_{|\vec n|\leq R}\frac{1}{|B|}\bigl(u(\vec n)- u(\vec n+\vec\ell)\bigr)^2&\leq&|\vec \ell|\sum_{\vec k\in\gamma_{\vec\ell}}\sum_{|\vec n|\leq R}\frac{1}{|B|}\sum_{i=1}^d\Bigl(\cr
&&\bigl(\partial_{i}^{*}u(\vec n+\vec k)\bigr)^2+\bigl(\partial_{i}u(\vec n+\vec k)\bigr)^2\Bigr).
\end{eqnarray*}
Or dans la somme intérieure en $\vec n$ du membre de droite, l'argument $\vec n+\vec k$ de $\partial u$
varie sur un ensemble contenu dans la boule de rayon $R+1$. Le membre de droite est donc majoré par
$$|\vec \ell|\sum_{\vec k\in\gamma_{\vec\ell}}\sum_{|\vec n|\leq R+1}\frac{1}{|B|}\sum_{i=1}^d\Bigl(\bigl(\partial_{i}^{*}u(\vec n)\bigr)^2
+\bigl(\partial_{i}u(\vec n)\bigr)^2\Bigr).
$$
La somme en $\vec n$ ne dépend plus 
 de $\vec k$.
Comme la somme en $\vec k$ porte sur $|\vec \ell|$ termes, c'est encore majoré par
$$|\vec \ell|^2\sum_{|\vec n|\leq R+1}\frac{1}{|B|}\sum_{i=1}^d|\partial_{i}u(\vec n)|^2+|\partial_{i}^*u(\vec n)|^2.
$$
Il reste à sommer en $\vec \ell$, qui est de norme $\leq 2R$. Comme cette somme porte sur $|B|$ termes,
cela démontre bien le lemme.\qed

Passons à la démonstration de la proposition~\ref{prharmo}. Rappelons la proposition~6.2 de~\cite{Delmotte}.
Bien qu'il n'y ait rien d'aléatoire dans l'article de Delmotte, nous conserverons notre notation
d'un milieu dépendant de $\omega$,
pour ne pas introduire de nouvelles notations. Notamment
 une fonction $u$ sera harmonique si elle vérifie ${\mathcal L}_\omega(u)=0$,
où ${\mathcal L}_\omega$ est l'opérateur introduit dans l'égalité~(\ref{ccharmo}).
\begin{prc}[Delmotte]
Il existe des constantes  $\alpha $ et $C_1$  ne dépendant que de la dimension $d$, et des bornes
d'ellipticité $a$ et $b$, telles que pour toute fonction $u$ harmonique définie sur ${\mathbb Z}^d$ on a:
$$\mbox{si }|\vec m|,\ |\vec n|<R \mbox{ alors }|u(\vec m)-u(\vec n)|\leq C_1\max_{|\vec \ell|\leq2R} |u(\vec\ell) |
\Bigl( \frac{|\vec m-\vec n|}{R}\Bigr)^\alpha.
$$
\end{prc}
Appliquée à la fonction $u$ définie sur ${\mathbb Z}^d$ par $u(\vec n)=S_{\vec n}(\omega)$, 
ou plus précisément à $u-\bar u_{4R}$, cette proposition s'écrit
$$\mbox{si }|\vec m|,\ |\vec n|<R \mbox{ alors }|S_{\vec m}(\omega)-S_{\vec n}(\omega)|\leq C(R,\omega)
\Bigl( \frac{|\vec m-\vec n|}{R}\Bigr)^\alpha.
$$
avec 
$$C(R,\omega)=C_1\max_{|\vec \ell|\leq2R} \Bigl|S_{\vec\ell}(\omega) -\frac{1}{|B_{4R}|}\sum_{\vec n\in B_{4R}}S_{\vec n}(\omega)
\Bigr|
.$$
Rappelons enfin la proposition 5.3 du même article.
\begin{prc}[Delmotte] Il existe une constante $C_2$  ne dépendant que de la dimension $d$, et des bornes
d'ellipticité $a$ et $b$, telle que si ${\mathcal L}_\omega u\geq 0$, alors
$$\max_{B_{2R}} u\leq C_2\Bigl(\frac{1}{|B_{4R}|}\sum_{\vec \ell\in B_{4R}}u(\vec\ell)^2\Bigr)^{\frac{1}{2}}.
$$
\end{prc}
D'après cette proposition, appliquée à la fonction $u-\bar u_{4R}$
 puis à son opposée $-u+\bar u_{4R}$,
on a
$$C(R,\omega)\leq C_1C_2\Bigl(\frac{1}{|B_{4R}|}\sum_{\vec \ell\in B_{4R}} \Bigl|S_{\vec\ell}(\omega) -\frac{1}{|B_{4R}|}\sum_{\vec n\in B_{4R}}S_{\vec n}(\omega)
\Bigr|
^2\Bigr)^{\frac{1}{2}}.
$$
 Enfin l'inégalité de Poincaré (lemme~\ref{poincarre}) assure que
ce majorant est lui-même  majoré par
$$
C'(R,\omega)=C_1C_2\Bigl(\frac{(4R)^2}{|B_{4R}|} 4\sum_{|\vec n|\leq 4R+1} \sum_{i=1}^d|f_i(T_{\vec n}\omega)|^2
+|f_i(T_{\vec n-\vec e_i}\omega)|^2
\Bigr)^{\frac{1}{2}}.
$$  
En divisant par $R$, on
 achève la démonstration de la proposition~\ref{prharmo} avec 
\begin{equation}\label{CC}
 C=\sup_{R\geq 1}\frac{1}{R}C'(R,\omega),
\end{equation}
%
 qui est presque sûrement finie d'après  le théorème ergodique ponctuel de Wiener
 sur les boules de ${\mathbb Z}^d$ (voir~\cite{Wiener}), et l'intégrabilité $L^2$ des fonctions $f_i$.
 \qed

Nous sommes maintenant en mesure de démontrer  le théorème~\ref{thharmo}.
Soit $(S_{\vec n})_{\vec n}$ le  cocycle donné par le théorème de Kunnemann.
Soit $\vec n\in{\mathbb Z}^d$.
 D'après l'équation~(\ref{cc1}) définissant les cocycles de degré~$1$, on a pour tout entier $k\geq 1$
 $$S_{k\vec n}=\sum_{j=0}^{k-1}S_{\vec n}\circ T_{\vec n}^j.
 $$
  D'après le théorème ergodique en dimension~$1$ 
 appliqué
  la fonction $S_{\vec n}$
  et la transformation $T_{\vec n}$, 
 la quantité 
 $$\frac{S_{k\vec n}-k\sum_{i=1}^dn_i\int_\Omega f_i\ dP}{k|\vec n|}
 $$
converge  donc presque sûrement, pour $k$ tendant vers l'infini.
La transformation $T_{\vec n}$ n'est pas nécessairement ergodique.
Mais d'après la remarque~\ref{re1} qui suit l'énoncé du théorème ergodique de Boivin-Derriennic, 
et l'hypothèse d'ergodicité de l'action $T$,
on a la convergence dans $L^1$ vers $0$.
Donc la limite ponctuelle vaut aussi $0$.

Soit un entier $n$ fixé. Le nombre de vecteur $\vec n$ à coordonnées entières vérifiant $|\vec n|=n$ est fini.
Donc d'après ce qui précède, avec probabilité $1$, pour tout $\varepsilon$, il existe un entier 
$R_0\geq 1$ 
tel que  
\begin{equation}\label{eps}
\mbox{ si }kn\geq R_0\mbox{ et si }|\vec n|=n\mbox{ alors }\Bigl|\frac{S_{k\vec n}(\omega)-
k\sum_{i=1}^dn_i\int_\Omega f_i\ dP
}{k n}\Bigr|\leq \varepsilon.
\end{equation}
 
 Notons $E_n$ l'ensemble des multiples des points de ${\mathbb Z}^d$ de norme $n$, soit
 $$E_n=\{k\vec n,\ k\in{\mathbb Z},\ |\vec n|=n\}.
 $$
On a alors le lemme élémentaire suivant
\begin{lemme}\label{l2}
Soit $n\geq 1$. Pour tout $\vec m\in{\mathbb Z}^d$ il existe $\vec \ell\in E_n$ tel que
$$|\vec\ell|\leq|\vec m| \mbox{ et } |\vec \ell-\vec m|\leq n+\frac{|\vec m|}{n}d.
$$
\end{lemme}
\par\noindent{\it Démonstration. --- } Considérons la division euclidienne de $|\vec m|$ par $n$
$$|\vec m|=kn+r,
$$
avec $0\leq r\leq n-1$. Le vecteur $\vec x=\frac{\vec m}{|\vec m|}n$ a pour norme $n$. Ses coordonnées
$(x_i)_{1\leq i\leq d}$ ne sont pas nécessairement
entières, mais il existe un vecteur $\vec n$ à coordonnées entières de norme $n$, à une distance $\leq d$ de
ce dernier. 
En effet, en notant $n'_i$ la partie entière de $x_i$, pour $i=1,\ldots,d$, et en posant
$$K=\sum_{i=1}^d(x_i-n'_i)=n-\sum_{i=1}^dn'_i,
$$
on voit que $K$ est un entier vérifiant $0\leq K\leq  d-1$. Considérons alors $n_i=n'_i+1$ si $i\leq K$ et $n_i=n'_i$ si $i\geq K+1$.
Le point $\vec n$ ainsi défini vérifie bien  $|\vec n|=n$ et $|\vec n-\vec x|\leq d$.

Posons $\vec \ell=k\vec n$. On a
$$|\vec m-\vec\ell|\leq \bigl|\vec m-\frac{\vec m}{|\vec m|}nk\bigr|+\bigl|\frac{\vec m}{|\vec m|}nk-\vec \ell\bigr|\leq r+kd,
$$
ce qui est bien majoré par $n+\frac{|\vec m|}{n}d$. Ceci démontre le lemme~\ref{l2}.\qed

Pour $0<\varepsilon<1/2 $ fixé, choisissons $n>d/\varepsilon$, 
puis $R_0>n/\varepsilon$. Avec ce choix, le lemme précédent  assure que pour tout $\vec m$ tel que $|\vec m|>R_0$
 il existe $\vec \ell\in E_n$ tel que $|\vec\ell|\leq|\vec m|$ et $|\vec \ell-\vec m|<2\varepsilon|\vec m|$.

On déduit alors de la proposition~\ref{prharmo} appliquée à $R=|\vec m|$ 
que
$$\bigl|\frac{S_{\vec m}(\omega)}{|\vec m|}-\frac{S_{\vec \ell}(\omega)}{|\vec m|}\bigr|\leq C
(2\varepsilon)^\alpha.
$$
De la majoration $\Bigl|\sum_{i=1}^d(m_i-\ell_i)\int_\Omega f_i\ dP\Bigr|\leq
|\vec m-\vec \ell|
\max_i|\int_\Omega f_i\ dP|
$ 
il découle que
$$\Bigl|\frac{S_{\vec m}(\omega)-\sum_{i=1}^dm_i\int_\Omega f_i\ dP}{|\vec m|}-\frac{S_{\vec \ell}(\omega)-\sum_{i=1}^d\ell_i\int_\Omega f_i\ dP}{|\vec m|}\Bigr|
$$
est majoré par
$$C
( 2\varepsilon)^\alpha+2\varepsilon\max_i\bigl|\int_\Omega f_i\ dP\bigr|.
$$
Enfin, on  peut choisir $R_0$ suffisamment grand pour que~(\ref{eps}) soit vérifié. Or, pour tout $\vec m$ tel que $|\vec m|>R_0\bigl(1-2\frac{d}{n}\bigr)^{-1}
$ et $\vec \ell$ donné par le lemme~\ref{l2},
on a $R_0\leq |\vec \ell|\leq |\vec m|$. Donc d'après~(\ref{eps}) et la ligne précédente
$$\Bigl|\frac{S_{\vec m}(\omega)-\sum_{i=1}^dm_i\int_\Omega f_i\ dP}{|\vec m|}\Bigr|\leq C
( 2\varepsilon)^\alpha+\Bigl(2\max_i\bigl|\int_\Omega f_i\ dP\bigr|+1\Bigr)\varepsilon.
$$
Comme $\varepsilon$ est arbitrairement petit, cela montre bien la convergence ponctuelle de 
$$\bigl|\frac{S_{\vec m}(\omega)-\sum_{i=1}^dm_i\int_\Omega f_i\ dP}{|\vec m|}\bigr|$$
 vers $0$, pour $|\vec m|$ tendant vers l'infini, ce qui démontre le 
théorème~\ref{thharmo}.\qed
\section{Question} La question naturelle à la suite de ce travail est: dans quelles conditions peut-on étendre
le théorème~\ref{thharmo} au cas non elliptique. 
Dans le cas où les conductances $c_i\circ T_{\vec n}$, pour $\vec n\in{\mathbb Z}^d$ et $i=1,\ldots,d$,
sont indépendantes, ce probléme est notamment abordé dans~\cite{BiskupPrescott}.
Dans le cas  où les conductances sont simplement supposées de loi stationnaire,
nous n'avons qu'une référence: l'article de 
Berger-Biskup~\cite{BergerBiskup}. Cet article se place dans le cadre de la percolation, c'est-à-dire 
avec des fonctions $c_i$ prenant  les valeurs $0$ ou $1$.  La convergence ponctuelle
du cocyle (souvent appelé \og correcteur\fg\ dans la méthode de la martingale évoquée ci-dessus) qui fait l'objet du théorème 5.1,
est démontré pour la dimension $2$. La méthode proposée ne passe pas à des dimensions plus grandes.

\bibliographystyle{smfalpha}
\providecommand{\bysame}{\leavevmode ---\ }
\providecommand{\og}{``}
\providecommand{\fg}{''}
\providecommand{\smfandname}{et}
\providecommand{\smfedsname}{\'eds.}
\providecommand{\smfedname}{\'ed.}
\providecommand{\smfmastersthesisname}{M\'emoire}
\providecommand{\smfphdthesisname}{Th\`ese}

\end{document}